\documentclass[11pt]{article}
\usepackage {amsmath}
\usepackage {amsfonts}
\usepackage {amsthm}
\usepackage {amssymb}
\usepackage {fullpage}
\usepackage {amscd}
\usepackage[all,cmtip]{xy}

\title{Heisenberg Idempotents on Unipotent Groups}
\author{Tanmay Deshpande\thanks{Partially supported by NSF grant DMS-0701106.}}
\date{}

\newtheorem {thm} {Theorem} [section]
\newtheorem {prop} [thm] {Proposition}
\newtheorem {lem} [thm] {Lemma}
\newtheorem {cor} [thm] {Corollary}
\theoremstyle{definition}
\newtheorem {defn} [thm] {Definition}
\theoremstyle{remark}
\newtheorem {rk} [thm]  {Remark}

\newcommand{\bpf}{\begin{proof}}
\newcommand{\epf}{\end{proof}}
\newcommand{\rar}[1]{\stackrel{#1}{\longrightarrow}}
\newcommand{\f}{\mathbb}

\renewcommand{\L}{\mathcal{L}}
\newcommand{\cM}{\L}
\newcommand{\ot}{\otimes}

\newcommand{\Q} {\mathbb{Q}}

\newcommand{\F} {\mathtt{k}}
\renewcommand{\k} {\mathtt{k}}
\newcommand{\K} {\mathbb{K}}
\newcommand{\M} {\mathcal{M}}
\newcommand{\talpha} {\tilde{\alpha}}
\newcommand{\D} {\mathcal{D}}
\newcommand{\C} {\mathcal{C}}
\newcommand{\U} {\tilde{U}}
\newcommand{\Hom} {\hbox{Hom}}

\newcommand{\Qlcl} {\overline{\mathbb{Q}}_l}
\newcommand{\beq}{\begin{equation}}
\newcommand{\eeq}{\end{equation}}
\newcommand{\bthm}{\begin {thm}}
\newcommand{\ethm}{\end {thm}}
\newcommand{\bprop}{\begin {prop}}
\newcommand{\eprop}{\end {prop}}
\newcommand{\bcor}{\begin {cor}}
\newcommand{\ecor}{\end {cor}}
\newcommand{\blem}{\begin{lem}}
\newcommand{\elem}{\end{lem}}
\newcommand{\bdefn}{\begin{defn}}
\newcommand{\edefn}{\end{defn}}
\newcommand{\brk}{\begin{rk}}
\newcommand{\erk}{\end{rk}}
\newcommand{\pt}{\hbox{pt}}

\newcommand{\noin}{\noindent}
\newcommand{\A} {\mathcal{A}}
\newcommand{\tgamma} {\tilde{\gamma}}
\newcommand{\tphi} {\tilde{\phi}}
\newcommand{\tg} {\tilde{\gamma}}

\newcommand{\g} {{\gamma}}

\newcommand{\ga} {{\gamma_1}}
\newcommand{\gb} {{\gamma_2}}
\newcommand{\cg} {{\ast}}

\newcommand{\HN}{H\ltimes N}
\newcommand{\V}{H/N}
\renewcommand{\l}{\Qlcl}
\newcommand{\id}{id}
\newcommand{\al}{\alpha}
\newcommand{\xrar}[1]{\xrightarrow{#1}}
\newcommand{\tmu}{\tilde{\mu}}
\newcommand{\Meg}{\M_e^{\Gamma}}
\newcommand{\bit}{\begin{itemize}}
\newcommand{\eit}{\end{itemize}}

\newcommand{\cG}{\mathcal{G}}
\newcommand{\cS}{\mathcal{S}}
\newcommand{\B}{\mathcal{B}}

\begin{document}
\maketitle

\section{Introduction} \label{introduction}
\noin 
Let $G$ be a possibly disconnected algebraic group over an algebraically closed field $\F$ of characteristic $p>0$,  such that its neutral connected component, $H=G^0$, is a unipotent group. We recall that an algebraic group over $\F$ is defined to be a smooth group scheme of finite type over $\F$. Let us fix a prime number $l\neq p$. If $X$ is a $\F$-scheme, we use $\D(X)$ to denote the bounded derived category of $\l$-complexes on $X$. If the group $G$ acts on $X$, we use $\D_G(X)$ to denote $G$-equivariant bounded derived category of $\l$-complexes on $X$.
 
\subsection{Heisenberg Idempotents}
Let $N$ be a closed connected normal subgroup of $G$, hence of $H$, such that the quotient $H/N$ is commutative. Let $\L$ be a $G$-equivariant multiplicative $\l$-local system on $N$. In particular, $\L$ is $H$-equivariant. For a $\F$-scheme $X$, let $X_{perf}$ denote its perfectization. Then as defined in \cite{B}, we get a map $\phi_\L:(\V)_{perf}\to (\V)_{perf}^*$, where $(\V)_{perf}^*$ is the Serre dual of $(\V)_{perf}$. We will only need to think about the $\F$-points of $(\V)^*_{perf}$ and these can be identified with multiplicative local systems on $\V$. Let $\L$ be such that the map $\phi_\L$ is an isogeny, i.e. such that $(N,\L)$ is an admissible pair for $G$, in the terminology of \cite{B}. Let $K_\L$ denote the kernel of this isogeny. Let $\D_G(G)$ denote the $G$-equivariant (under conjugation action) bounded derived category of $\l$-complexes on $G$ and $\D_H(G)$ the $H$-equivariant bounded derived category. The categories $\D(G), \D_H(G)$ and $\D_G(G)$ have the structure of a monoidal category under convolution of complexes and $\D_G(G)$ has the structure of a braided monoidal category. As described in \cite{B}, in this situation, we can define a closed idempotent $e\in \D_G(G)$. More explicitly, $e= \L\otimes \f{K}_N$ considered as a complex on $G$ by extending by zero outside $N$, with the $G$-equivariant structure coming from the $G$-equivariant structure on $\L$. Here $\f{K}_N=\Qlcl[2\dim N](\dim N)$ is the dualizing complex on $N$. An idempotent on $G$ obtained in this way is known as a Heisenberg idempotent. In this situation, we would like to study the Hecke subcategory $e\D_G(G)$. 

\subsection{Main Results}

In this article, we will describe the category $e\D_G(G)$. First, we work with the category $e\D_H(G)$. Since $H$ is connected, the forgetful functor from $\D_H(G)$ to $\D(G)$ is fully faithful. Hence we often implicitly consider the categories $\D_H(G)$ and $e\D_H(G)$ as full subcategories of $\D(G)$. Note that since $e\in \D_G(G)$, we have $eM\cong Me$ for all $M\in \D(G)$. Moreover, $e$ is a closed idempotent. Hence, it follows from \cite[\S2]{BD}  that $e\D_H(G)=e\D_H(G)e$ is a monoidal category with unit object $e$. Let $\M^{{perv}}_e$ denote the full subcategory of $e\D_H(G)$ consisting of perverse sheaves. Then $\M^{{perv}}_e$ is an additive $\Qlcl$-linear subcategory of $\D_H(G)$. It is also clear that $e\in \M^{{perv}}_e[\dim N]$. We will prove the following results conjectured by V. Drinfeld.
\bthm\label{Main1}
$\M^{{perv}}_e$ is a semisimple abelian category with finitely many simple objects (up to isomorphism). Each simple object in $\M^{perv}_e$ is a suitably shifted local system supported on a closed non-singular subvariety of $G$. Moreover, the canonical functor from $D^b(\M^{{perv}}_e)$ to $e\D_H(G)$ is an equivalence.
\ethm
\noin We prove this result in \S\ref{secmain1}. To prove the theorem, it will be convenient to use the notion of quasi-equivariant complexes, which we define in \S\ref{thecatdulx}. In \S\ref{supportofqe}, we describe the support of such quasi-equivariant complexes. In \S\ref{qecom}, we describe the category of quasi-equivariant complexes on a homogeneous space. In \S\ref{altdes}, we give an alternative description of the category $e\D_H(G)$ as the category of certain quasi-equivariant complexes on $G$ with respect to the action of $H\ltimes N$ on $G$, where $H$ acts by conjugation and $N$ acts by left multiplication. In \S\ref{supportedhg}, we show that these quasi-equivariant complexes can only be supported on  finitely many $\HN$-orbits in $G$. Each orbit is a non-singular closed subvariety of $G$ and is a homogeneous space for $\HN$. Hence, the results from \S\ref{qecom} will give us an explicit description of the category of quasi-equivariant complexes supported on a single orbit. In particular, we will show that the category of quasi-equivariant perverse sheaves supported on a single orbit is semisimple abelian with finitely many simple objects and that the category of quasi-equivariant complexes supported on that orbit is the bounded derived category of the former category. (See Proposition \ref{transact}.)
\bthm\label{Main2}
The subcategory $\M_e:=\M^{{perv}}_e[\dim N]\subset \D_H(G)$ is closed under convolution, and is a monoidal category with unit object $e$.
\ethm
\noin Since $e$ is a unit object in $e\D_H(G)$, the theorem is equivalent to the assertion that if $M,L\in\mathcal{M}_e^{perv}$, then $M*L\in\mathcal{M}_e^{perv}[-\dim N]$. We prove this theorem in \S\ref{potm2}. Let us explain the idea of the proof. In \S\ref{acoat}, we show using Artin's theorem, that $M\ast L\in {}^p\D^{\geq \dim N}(G)$. Then, we will need to use a notion of duality in the category $e\D_H(G)$, which is weaker than rigidity. We describe this in \S\ref{diedhg}. Namely,  we construct an antiequivalence $L\mapsto L^{\vee}$ of $e\D_H(G)$, such that for $M,L\in e\D_H(G)$, we have functorial isomorphisms $$\Hom(M,L^{\vee})=\Hom(M\ast L, e).$$
The subcategory $\M_e$ is stable under this antiequivalence.
In \S\ref{cwtd}, we  compute $M^{\vee}\ast M$, where $M\in \M_e$ is a simple object, and see that it lies in $\M_e$. Moreover, we will see that it is supported on $H$. In \S\ref{tcedhh}, we describe all the simple objects of the category $\M_e$ supported on $H$ (Proposition \ref{edhh}) and compute the convolution of a general complex in $\M_e$ with these simple objects. In particular, we will see that these convolutions also lie in $\M_e$. (See Proposition \ref{edhhmod}.) From this, it follows that for any $M,L\in \M_e$, $M^{\vee}\ast(M\ast L)\cong(M^{\vee}\ast M)\ast L\in \M_e$. Then in \S\ref{potm2}, using the semisimplicity of $\M_e$, we will deduce that we must in fact have that $M\ast L\in \M_e$.

\bthm\label{Main3}
The category $\M_e$ is rigid monoidal, and hence a fusion category.
\ethm
\noin Let $\Gamma=G/H$. Then we have the $\Gamma$-grading $\D(G)=\bigoplus\limits_{H\tgamma \in G/H}\D(H\tgamma)$. This gives us the grading
$\M_e=\bigoplus\limits_{\gamma \in \Gamma}\M_{e,\gamma}$, where $\M_{e,\gamma}$ is the full subcategory of $\M_e$ consisting of objects supported on the $H$-coset corresponding to $\gamma\in \Gamma$. We will see that the trivial component $\M_{e,1}$ is pointed, i.e. all simple objects in $\M_{e,1}$ have an inverse.
In \S\ref{rigidity}, we prove that under these conditions, a tensor category satisfying the weak duality property described above is in fact rigid. From this we can conclude that $\M_e$ is a fusion category.

We will see in \S\ref{brcross} that the categories $\D_H(G)$ and $e\D_H(G)$ have the structure of a braided $\Gamma$-crossed category. This induces a braided $\Gamma$-crossed structure on $\M_e$. The $\Gamma$-equivariantizations $e\D_H(G)^\Gamma$ and $\Meg$ are braided monoidal categories. It is easy to see that we have the following:
\blem
We have an equivalence $e\D_G(G)\cong e\D_H(G)^{\Gamma}$of braided monoidal categories. Under this equivalence, the full subcategory of $e\D_G(G)$ consisting of those objects whose underlying $\Qlcl$-complex is a perverse sheaf shifted by $\dim N$, gets identified with $\Meg$. 
\elem

As defined in \cite{B}, we have the twist automorphism $\theta$ of the identity functor on $\D_G(G)$. This gives us a twist, which we also denote by $\theta$, in  $\Meg$.
\bthm\label{Main4}
\bit 
\item[(i)]The category $\Meg$ is a semisimple abelian category with finitely many simple objects. For each simple object, the underlying $\Qlcl$-complex is a suitably shifted local system supported on a closed non-singular subvariety of $G$. The canonical functor $D^b(\Meg)\to e\D_H(G)^\Gamma\cong e\D_G(G)$ is an equivalence.
\item[(ii)]The category $\Meg$ is a braided fusion category.
\item[(iii)]The twist $\theta$ defines a ribbon structure on $\Meg$ and hence gives $\Meg$ the structure of a pre-modular category. In fact, $\Meg$ is a modular category.
\eit
\ethm
\noin Statements (i) and (ii) follow readily from the previous results. We verify in \S\ref{ribbon} that $\theta$ indeed defines a ribbon structure on $\Meg$. In \S\ref{titcedhh}, we show that the twists in the category $e\D_H(H)$ define a quadratic form (which we also denote by $\theta$) $\theta:K_\L\to \Qlcl^*$ which gives us a polarization of a certain non-degenerate symmetric pairing $B:K_\L\times K_\L\to \Qlcl^*$. From this we conclude that the category $\M_{e,1}$ is the modular category corresponding to the quadratic form $\theta$ on $K_\L$. In particular, $\M_{e,1}$ is a non-degenerate braided fusion category. Using results from \cite[\S4.4.8]{DGNO}, we can deduce that $\Meg$ is a non-degenerate braided fusion category, and hence a modular category.

\subsection{Acknowledgments}
I would like to thank my advisor, V. Drinfeld and M. Boyarchenko for introducing me to this subject, for the many useful discussions and for suggesting various improvements and corrections.

\section{The categories $\D(G), \D_G(G)$ and $\D_H(G)$} \label{catdg}
\subsection{Convolution of complexes} \label{convolution}
Let $\mu:G\times G\to G$ denote the group operation on $G$. For any algebraic group $G$, we have a convolution with compact supports, which is a bifunctor $\cg:\D(G)\times \D(G)\to \D(G)$. Namely, for $M,L\in \D(G)$, we set $M\cg L = \mu_!(M\boxtimes L)$. This makes each of $\D(G), \D_G(G)$ and $\D_H(G)$ a monoidal category with the unit object ${\mathbf 1}$ given by the delta sheaf supported at the identity $1$ of $G$, with the stalk at $1$ equal to $\Qlcl$. If $L$ is a $G$-equivariant complex, then we can define isomorphisms $M\cg L\cong L\cg M$ for all $M\in \D(G)$. In particular, we have $e\cg M\cong M\cg e$ for all $M\in \D(G)$.

For $M\in \D(G)$ and $g\in G$, let us denote by $M^g$ the right translate of $M$ by $g$, i.e. $M^g=r_{g^{-1}}^*M$, where $r_{g^{-1}}:G\to G$ denotes multiplication on the right by $g^{-1}$. Similarly, we define ${}^gM=l_{g^{-1}}^*M$, where $l_{g^{-1}}$ denotes left multiplication by $g^{-1}$. It is easy to check that $^gM\cong\delta_g*M$ and $M^g\cong M*\delta_g$, where $\delta_g$ denotes the delta-sheaf at $g$. This observation implies the following:
\bprop\label{properties}
Let $M,L\in \D(G)$ and $g, g_1,g_2\in G$. Then we have the following canonical isomorphisms:
\begin{itemize}
\item[(a)] $(M^{g_1})^{g_2}\cong M^{g_1g_2}$ and ${}^{g_1}(^{g_2}M)\cong {}^{g_1g_2}M$.
\item[(b)] $( {}^{g_1}M)^{g_2}\cong {}^{g_1}(M^{g_2})$.
\item[(c)] $M\ast (L^{g})\cong (M\ast L)^g$ and $({}^gM)\ast L\cong {}^g(M\ast L)$.
\end{itemize}
\eprop

\subsection{Duality in $\D(G)$ and $\D_H(G)$}
Let us now describe a duality in $\D(G)$ and $\D_H(G)$, which is weaker than the notion of rigidity and rigid duals. Namely, there exists an antiequivalence $\f{D}_G^-={\f{D}}^-:\D(G)\to \D(G)$ such that we have functorial isomorphisms 
\beq\label{duality}
\Hom(M, \f{D}^-L)=\Hom(M\cg L, \mathbf{1}), \hbox{ for } M, L \in \D(G).
\eeq
This also induces an antiequivalence $\f{D}_G^-={\f{D}}^-:\D_H(G)\to \D_H(G)$ satisfying the above property for $M,L\in \D_H(G)$.
Explicitly, $\f{D}^- = \f{D}\circ \iota^* = \iota^*\circ \f{D}$, where $\iota:G\to G$ is the inversion map and $\f{D}:\D(G)\to \D(G)$ is the Verdier duality functor. 
\brk\label{verdierperverse}
Since the Verdier dual of a perverse sheaf is perverse, we see that the antiequivalence $\f{D}^-$ stabilizes the full subcategories of $\D(G), \D_H(G)$ consisting of perverse sheaves. 
\erk

\subsection{Duality in $e\D_H(G)$} \label{diedhg}
\noin Using the duality in $\D_H(G)$ described above, we can define a similar duality in the Hecke subcategory $e\D_H(G)$.

\bprop
The full subcategory $e\D_H(G)\subset \D_H(G)$ is stable under the functor $\f{D}^-$.
\eprop
\bpf
Suppose $L\in e\D_H(G)$ i.e. the map $L\to e\cg L$ is an isomorphism. Then from (\ref{duality}), we conclude that for all $M\in \D_H(G)$, the map $M\to M\cg e$ induces a bijection, $\Hom(M\cg e, \f{D}^{-}L)\cong \Hom(M, \f{D}^-L)$. From this we conclude that $\f{D}^-L\in \D_H(G)e=e\D_H(G)$. (See \cite{BD},\S2.)
\epf

\bprop\label{dualityinhecke}
We have an antiequivalence, $L\mapsto L^{\vee}$ of $e\D_H(G)$ such that we have functorial isomorphisms
\beq
\Hom(M, L^{\vee})=\Hom(M\cg L, e), \hbox{ for } M, L \in e\D_H(G).
\eeq
Namely, for any $L\in e\D_H(G)$, set $L^{\vee}=(\f{D}^{-}L)[2\dim N](\dim N)$. 
\eprop
\bpf
We first note that we have a canonical isomorphism $\f{D}^{-}e\cong e[-2\dim N](-\dim N)$, and hence a canonical isomorphism $e^{\vee}\cong e$. By (\ref{duality}), we have that $$\Hom(M, \f{D}^-L)=\Hom(M\cg L, \mathbf{1})=\Hom(M\ast L\ast e, \mathbf{1})=\Hom(M\ast L, \f{D}^{-}e).$$ Hence we conclude that we must have functorial isomorphisms $\Hom(M, L^{\vee})=\Hom(M\cg L, e)$. 
\epf

\brk\label{mestab}
Using the fact that $\f{D}^-$ stabilizes the full subcategory of perverse sheaves, we see that the antiequivalence ${(\cdot)}^{\vee}$ stabilizes $\M_e^{perv}[\dim N]=\M_e$. We also note that if $L\in e\D_H(G)$ is supported on $Z\subset G$, then $L^{\vee}$ is supported on  $\iota(Z)$.
\erk

\subsection{Braided $\Gamma$-crossed structure on $\D_H(G)$} \label{brcross}
Let us first recall some definitions. 
\bdefn\label{grading}
Let $\C$ be an additive monoidal category and let $\Gamma$ be a finite group. A $\Gamma$-grading on $\C$ is a decomposition $\C=\bigoplus\limits_{\gamma\in \Gamma}\C_\gamma$ such that for $\ga,\gb\in \Gamma$ we have $\C_{\ga}\otimes\C_{\gb}\subset \C_{\ga\gb}$. We say that a grading is faithful if $\C_\gamma\neq 0$ for all $\gamma\in \Gamma$.
\edefn

\bdefn\label{braidedcrossed}
Let $\Gamma$ be a finite group. A braided $\Gamma$-crossed category $\C$ is an additive monoidal category $\C$ equipped with the following structures:
\bit
\item[(i)]a $\Gamma$-grading $\C=\bigoplus\limits_{\gamma\in \Gamma}\C_\gamma$;
\item[(ii)]a monoidal action of $\Gamma$ on $\C$ such that $g(\C_h)\subset \C_{ghg^{-1}}$ for all $g,h\in \Gamma$;
\item[(iii)]for $g\in \Gamma$, $X\in \C_g, Y\in \C$ isomorphisms \beq c_{X,Y}:X\otimes Y\to g(Y)\otimes X\eeq functorial in $X$ and $Y$ called $\Gamma$-braiding isomorphisms satisfying the following conditions:
\bit
\item[(a)]$\gamma(c_{X,Y})=c_{\g(X),\g(Y)}$ for all $\g\in \Gamma$;
\item[(b)]the following diagrams commute for all $g,h\in \Gamma$, $X\in \C_g$ and $Y\in \C_h$
\begin{equation} 
\xymatrix{X\ot Y\ot Z\ar[rr]^{c_{X,Y\ot Z}} \ar[d]_{c_{X,Y}\ot \id_Z
}&&g(Y\ot Z)\ot X\ar[d]^{\cong}\\ g(Y)\ot X\ot Z\ar[rr]^{\id_{g(Y)}\ot c_{X,Z}}&&g(Y)\ot g(Z)\ot X}
\end{equation}
\beq
\xymatrix{X\ot Y\ot Z\ar[rr]^{c_{X\ot Y,Z}} \ar[rd]_{\id_X\ot c_{Y,Z}
}&&gh(Z)\ot X\ot Y\\ &X\ot h(Z)\ot Y .\ar[ru]_{c_{X,h(Z)}\ot \id_Y}&}
\eeq
\eit
\eit
\edefn

Let us now describe the braided $\Gamma$-crossed structure on $\D_H(G)$. Firstly, we have the grading $\D_H(G)=\bigoplus\limits_{H\tgamma \in G/H}\D_H(H\tgamma)$. Then we also have that $$(\tg^{-1})^*(\D_H(H\tg'))\subset \D_H(H\tg\tg'\tg^{-1})$$ for $\tg,\tg'\in G$, where $\tg^{-1}:H\tg\tg'\tg^{-1} \to H\tg'$ denotes conjugation by $\tg^{-1}$. For $\g\in \Gamma$, let $\tg\in G$ denote a lift. The functors $(\tg^{-1})^*:\D_H(G)\to \D_H(G)$ for $\g\in \Gamma$ define an action of the finite group $\Gamma$ on $\D_H(G)$.

Let us now construct the crossed braiding, namely for $M\in \D_H(H\tg)$ and $L\in \D_H(G)$ we construct functorial isomorphisms 
\beq
c_{M,L}:M\ast L \stackrel{\cong}{\longrightarrow} (\tg^{-1})^*L\ast M.
\eeq
Note that we have the following commutative diagram:

$$\xymatrix{H\tg\times G\ar[d]_{\tau} \ar[r]^{\xi} &H\tg\times G\ar[d]^{\mu}\\
G\times H\tg\ar[r]_{\mu} &G ,}$$
where $\tau(h\tg, g)=(g, h\tg)$ and $\xi(h\tg, g)= (h\tg, \tg^{-1}h^{-1}gh\tg)$. Hence we have 
$$(\tg^{-1})^*L\ast M\cong \mu_!\left((\tg^{-1})^*L\boxtimes M\right)$$
$$\cong \mu_!\tau_!\left(M\boxtimes(\tg^{-1})^*L \right)$$
$$\cong \mu_!\xi_!\left(M\boxtimes(\tg^{-1})^*L \right).$$ Hence it is enough to construct an isomorphism $\xi_!\left(M\boxtimes(\tg^{-1})^*L \right)\stackrel{\cong}{\longrightarrow}M\boxtimes L$ of complexes on $H\tg\times G$, or in other words, an isomorphism $p_1^*M\otimes p_2^*(\tg^{-1})^*L\stackrel{\cong}{\longrightarrow} \xi^*p_1^*M\otimes \xi^*p_2^*L\cong p_1^*M\otimes C'^*(\tg^{-1})^*L$, where $p_1, p_2$ denote the two projections from $H\tg\times G$ and $C':H\tg\times G\to G$ is defined by $(h\tg, g)\mapsto h^{-1}gh$. Such an isomorphism is defined using the $H$-equivariant structure on $(\tg^{-1})^*L$. Hence we get a braided $\Gamma$-crossed structure on $\D_H(G)$. This also defines a braided $\Gamma$-crossed structure on $e\D_H(G)$.

\section{Quasi-equivariant complexes} \label{qeqcom}
In this section, we will describe the notion of quasi-equivariant complexes and give descriptions of $e\D(G)$ and $e\D_H(G)$ as categories of certain quasi-equivariant complexes on $G$.

\subsection{Multiplicative local systems}
Let us first define the notion of a multiplicative local system on a possibly disconnected algebraic group.
\bdefn\label{d:multiplicative-local-system}
Let $G$ be an algebraic group over $\k$ and let $\mu:G\times
G\to G$ denote the multiplication morphism. A
multiplicative local system on $G$ is a pair $(\L,\beta)$,
where $\L$ is a nonzero $\Qlcl$-local system on $G$ and
$\beta:\mu^*\L\stackrel{\cong}{\longrightarrow}\L\boxtimes\L$
is an isomorphism such that the two induced isomorphisms
$(\mu\times id_G)^*\mu^*\L\rar{\cong}\cM\boxtimes\cM\boxtimes\cM$
and
$(id_G\times\mu)^*\mu^*\cM\rar{\cong}\cM\boxtimes\cM\boxtimes\cM$
are equal modulo the canonical identification
$(\mu\times id_G)^*\mu^*\cM\cong (id_G\times\mu)^*\mu^*\cM$.
We will often abuse notation and only use $\L$ to denote a multiplicative local system.
\edefn

\brk
If $(\L,\beta)$ is a multiplicative local system on $G$, then $\beta$ induces an
isomorphism between the stalk of $\L$ at 1 and the 1-dimensional space
$\overline{\mathbb{Q}}_l$. Moreover, if $G$ is connected, then a rank
1 local system $\L$ on $G$ has a multiplicative structure if and only if
$\mu^*\L\cong \L\boxtimes \L$, and in this case, multiplicative
structures on $\L$ are in bijection with trivializations of the stalk
$\L_1$. Hence if the group $G$ is connected, we will not explicitly mention the multiplicative structure.
\erk

\brk\label{multloc}
Let $U$ be a unipotent algebraic group over $\k$. Let us fix an embedding $\psi:\Q_p/\f{Z}_p\to {\Qlcl}^*$. Then we can identify the group of isomorphism classes of central extensions of $U$ by $\Q_p/\f{Z}_p$ with the group of isomorphism classes of multiplicative local systems on $U$. (See \cite[\S5]{B}.) In particular, every multiplicative local system $\L'$ on $U$ comes from a central extension $0\to A\to \U\to U\to 0$ of $U$ by a finite group $A$, and a character $\chi:A\to \Qlcl^*$.
\erk
\medbreak

\subsection{The category $\D_{U,\L'}(X)$} \label{thecatdulx}
Suppose we have a unipotent group $U$ over $\k$ acting on a variety $X$ over $\k$. Let $(\L',\beta)$ be a multiplicative local system on $U$. Let us now define the category of $(U,\L')$-equivariant complexes on $X$.
\bdefn
Let $U, \L'$ and $X$ be as above. Let $\alpha:U\times X\to X$ denote the action.  By $\D_{U,\L'}(X)$, we denote the category of $(U,\L')$-equivariant complexes on $X$, whose objects are pairs $(M,\phi)$, where $M\in \D(X)$ and $\phi:\alpha^*M\stackrel{\cong}{\longrightarrow} \L'\boxtimes M$ such that the composition of the isomorphisms
\[
(\id_U\times\al)^*\al^* M \cong (\mu\times\id_X)^*\al^* M \xrar{\
(\mu\times\id_X)^*(\phi)\ } (\mu\times\id_X)^*(\cM'\boxtimes M)
\]
\[
\xrar{\ \beta\boxtimes\id_M\ } \cM'\boxtimes\cM'\boxtimes M
\]
equals the composition
\[
(\id_U\times\al)^*\al^* M \xrar{\ (\id_U\times\al)^*(\phi)\ }
(\id_U\times\al)^*(\cM'\boxtimes M) \cong \cM'\boxtimes(\al^* M)
\]
\[
\xrar{\ \id_{\cM'}\boxtimes\phi\ } \cM'\boxtimes\cM'\boxtimes M.
\]  A morphism $\nu:(M,\phi)\to (L, \psi)$ is a morphism $\nu:M\to L$ such that $\psi\circ\alpha^*(\nu)=(id_{\L'}\boxtimes\nu)\circ\phi$. The composition of two morphisms in $\D_{U,\L'}(X)$ is defined as their composition in $\D(X)$. Let us denote the category of $(U,\L')$-equivariant perverse sheaves on $X$ by $\M^{perv}_{U,\L'}(X)$.
\edefn

\brk\label{rk}
Note that $\D_{U,\Qlcl}(X)$ is exactly the category $\D_U(X)$ of $U$-equivariant complexes on $X$. Note that if $U$ is not unipotent, the above definition becomes unreasonable already in the case when $\L'$ is trivial.
\erk

\brk\label{connected}
We have a natural forgetful functor from $\D_{U,\L'}(X)$ to $\D(X)$. If $U$ is connected, this functor is fully faithful and its essential image is the full subcategory of $\D(X)$ consisting of complexes $M$ such that $\alpha^*M\cong \L'\boxtimes M$.
\erk


\bprop\label{centralext}
Let $U, (\L',\beta), X$ be as above. Consider a central extension $$0\to A \rar{i} \U\rar{\pi} U\to 0$$ of $U$ by a finite commutative group $A$.  Let $\tilde{\alpha}=\alpha\circ(\pi\times id_X):\U\times X\to X$ be the induced action. Let $\tilde{\mu}$ denote the multiplication on $\U$. Then we have a natural fully faithful functor $\D_{U,\L'}(X)\to \D_{\U,\pi^*\L'}(X)$ defined by $(M,\phi)\mapsto (M,(\pi\times id_X)^*\phi)$. The essential image of this functor consists of the objects $(M,\tilde{\phi})\in \D_{\U,\pi^*\L'}(X)$ such that $\tilde{\phi}|_{A\times X}:\talpha^*M|_{A\times X}\cong{\Qlcl}_A\boxtimes M\to {\Qlcl}_A\boxtimes M\cong \pi^*\L'|_A\boxtimes M$ is the identity.
\eprop
\bpf
Note that $(\pi^*\L',(\pi\times \pi)^*\beta)$ is a multiplicative local system on $\U$. It is clear that we indeed have such a functor. That this functor is fully faithful follows immediately from the definition of morphisms between quasi-equivariant complexes and the fact that two morphisms are equal in $\D(U\times X)$ if and only if their pullbacks via the \'{e}tale cover $(\pi\times id_X)$ are equal in $\D(\U\times X)$. Note that we have a canonical trivialization $\pi^*\L'|_{A}\cong {\Qlcl}|_A$ as multiplicative local systems. Let $(M,\phi)\in \D_{U,\L'}(X)$. Then $((\pi\times id_X)^*\phi)|_{A\times X}=(i\times id_X)^*(\pi\times id_X)^*\phi=(1\times id_X)^*\phi$ must be the identity. On the other hand, let $(M,\tilde{\phi})\in \D_{\U,\pi^*\L'}(X)$ be such that $\tphi|_{A\times X}$ is the identity. Now $\tphi$ satisfies the compatibility relation $$((\pi\times \pi)^*\beta\boxtimes id_M)\circ (\tmu\times id_X)^*\tphi=(id_{\pi^*\L'}\boxtimes \tphi)\circ (id_{\U}\times \talpha)^*\tphi.$$ Restricting this equality to $A\times \U\times X$, we deduce that $\tphi$ is a morphism in $\D_A(\U\times X)$ and hence $\tphi=(\pi\times id_X)^*\phi$ for some isomorphism $\phi:\alpha^*M\to \L'\boxtimes M$. Since $\tphi$ satisfies the compatibility relation, it follows that $\phi$ must also do so.
\epf

\subsection{Support of quasi-equivariant complexes} \label{supportofqe}

\bprop\label{support}
Let $(M,\phi)$ be a $(U,\L')$-equivariant complex on $X$. Let $x\in X$ be such that $M_x\neq 0$. Let $U_x$ denote the stabilizer of $x$ in $U$. Then we must have $\L'|_{U_x}\cong \Qlcl$.
\eprop
\bpf
We have an isomorphism $\phi:\alpha^*M\stackrel{\cong}{\longrightarrow} \L'\boxtimes M$ of complexes on $U\times X$. Restricting this isomorphism to $U_x\times \{x\}$, we get an isomorphism $M_x\cong \L'|_{U_x}\otimes M_x$ of complexes on $U_x\times \{x\}$, where we use $M_x$ to denote the constant complex. Since $M_x\neq 0$, we conclude that we must have $\L'|_{U_x}\cong \Qlcl$.
\epf

\subsection{Quasi-equivariant complexes on a homogeneous space} \label{qecom}
Following a suggestion made by M. Boyarchenko (Proposition \ref{transact}), let us now describe the category of quasi-equivariant complexes on a homogeneous space for a unipotent group $U$. Let $Vec = Vec_{\Qlcl}$ denote the category of finite dimensional vector spaces over $\Qlcl$ and $D^b(Vec)$ its bounded derived category. If $A$ is a finite group, the category $\D_A(\hbox{pt})$ is equivalent to the category whose objects consist of objects of $D^b(Vec)$ with an action of $A$. Here $\hbox{pt}$ stands for $\hbox{Spec}(\k)$ equipped with the trivial action by $A$.

\blem\label{lemmatransact}
Let $U$ be a unipotent group acting transitively on a variety $X$. Let $x\in X$. Let $U_x\subset U$ be the stabilizer of $x$. Then taking the stalk at $x$ induces an equivalence of categories $\D_U(X)\cong \D_{\pi_0(U_x)}(\hbox{pt})$. Under this equivalence, $\M^{perv}_{U}(X)[-\dim X]$ corresponds to $Rep(\pi_0(U_x))$.
\elem
\bpf
We will use Lemma 4.4 from \cite{B}. Consider the action of $U\times U_x$ on $U$, given by $(g,h)\cdot u = guh^{-1}$.  Let $N_1 = U\times \{1\}$ and $N_2 = \{1\}\times U_x$. The map $U\to \{x\}$ is an $N_1$-torsor, so $D_{U\times U_x}(U) \cong D_{U_x}(x)$ by Lemma 4.4 in \cite{B}. Note that $N_1$ admits a complement $H=\{(h,h)|h\in U_x\}\subset U\times U_x$. The map $\sigma:\{x\}\to U$ which sends $x$ to 1 is an $H$-equivariant section. Hence the equivalence above is induced by $\sigma^*$(by Lemma 4.5 from \cite{B}).
Also, the map $U\to X$ given by $g \mapsto gx$, is an $N_2$-torsor, so we have the quivalence $D_{U\times U_x}(U) \cong D_U(X)$ induced by pullback along the torsor map. 
So we see that we have an equivalence $D_U(X) \cong D_{U_x}(x)$ induced by taking the stalk at $x$. Now an object of $D_{U_x}(x)$ is an object of $D^b(Vec)$ with an action of $\pi_0(U_x)$. Moreover, since the action of $U$ on $X$ is transitive, $\M^{perv}_{U}(X)[-\dim X]$ is the full subcategory consisting of all equivariant local systems. Hence this subcategory corresponds to $Rep(\pi_0(U_x))$.
\epf

Let $U$ be a unipotent group over $\k$ and let $(\L',\beta)$ be a multiplicative local system on $U$. By Remark \ref{multloc}, there exists a central extension $0\to A\rar{i} \U\rar{\pi} U\to 0$ of $U$ by a finite commutative group $A$, along with an isomorphism $\tau:\pi^*\L'\rar{\cong} {\Qlcl}_{\U}$ of multiplicative local systems on $\U$. On the other hand, since $\pi\circ i=1$, we have a natural trivialization $\pi^*\L'|_A\cong {\Qlcl}_A$ of multiplicative local systems on $A$.  Hence we have an automorphism ${\Qlcl}_A\cong \pi^*\L'|_A \rar{\tau|_A} {\Qlcl}_A$ of the trivial multiplicative local system on $A$, or in other words, a homomorphism $\chi:A\to \Qlcl^*$. Suppose $U$ (and hence $\U$) acts transitively on $X$. For $x\in X$, let $U_x$  (respectively $\U_x$) be the stabilizer of $x$ in $U$ (respectively $\U$), so that we have a central extension $0\to A\to \U_x\to U_x\to 0$.

\bprop\label{transact}
Using the terminology of the paragraph above, we have an equivalence of categories $\D_{U,\L'}(X)\cong \D_{\pi_0(\U_x)}^{\chi}(\pt)$, where $D_{\pi_0(\U_x)}^{\chi}(\pt)\subset D_{\pi_0(\U_x)}(\pt)$ is the full subcategory consisting of objects of $D_{\pi_0(\U_x)}(\pt)$ such that $A$  acts\footnote{Note that we have a homomorphism from $A$ to $\pi_0(\U_x)$. Moreover, this homomorphism is injective if $\L'|_{U_x^0}$ is trivial.} by the character $\chi$. Under this equivalence, $\M^{perv}_{U,\L'}(X)[-\dim X]$ corresponds to the category $Rep_\chi(\pi_0(\U_x))$ of $\l$-representations of $\pi_0(\U_x)$ such that $A$ acts by the character $\chi$. Hence the canonical functor from $D^b(\M^{perv}_{U,\L'}(X))$ to $\D_{U,\L'}(X)$ is an equivalence. All the objects of $\M^{perv}_{U,\L'}(X)$ are local systems shifted by $\dim X$.
\eprop
\bpf
Note that we have a sequence of functors $\D_{U,\L'}(X)\to \D_{\U,\pi^*\L'}(X)\rar{\cong}\D_{\U}(X)\to \D_A(X)$. Let $(M,\phi)\in \D_{U,\L'}(X)$. Then the sequence of functor sends $$(M,\phi)\mapsto (M,(\pi\times id_X)^*\phi)\mapsto (M,(\tau\boxtimes id_M)\circ (\pi\times id_X)^*\phi)$$ 
$$\mapsto \left(M,(\tau|_A\boxtimes id_M)\circ (i\times id_X)^*(\pi\times id_X)^*\phi\right)$$
$$=\left(M,(\tau|_A\boxtimes id_M)\circ (1\times id_X)^*\phi\right).$$
The last object lies in $\D_A(X)$ and is given by an object $M\in \D(X)$ and the isomorphism $${\Qlcl}_A\boxtimes M\xrar{(1\times id_X)^*\phi} \pi^*\L'|_A\boxtimes M\xrar{\tau|_A\boxtimes id_M}{\Qlcl}_A\boxtimes M.$$ Note that the first map above comes from the natural trivialization of $\pi^*\L'|_A$, hence this last object (which lies in $\D_A(X)$) corresponds to the action of $A$ on the object $M\in \D(X)$ by the character $\chi:A\to \Qlcl^*$. By Proposition \ref{centralext}, the functor from $\D_{U,\L'}(X)$ to $\D_{\U,\pi^*\L'}(X)$ is fully faithful.  Hence we have a fully faithful functor from $\D_{U,\L'}(X)$ to $\D_{\U}(X)$. As we have seen above, the essential image of this functor is contained in the full subcategory of $\D_{\U}(X)$ consisting of objects on which $A$ acts by $\chi$. Moreover, it follows from Proposition \ref{centralext} that the essential image is precisely this full subcategory. Finally, using Lemma \ref{lemmatransact}, we can identify this full subcategory with $D_{\pi_0(\U_x)}^{\chi}(\pt)$. The remaining statements in the proposition are also clear.
\epf

\subsection{Alternative descriptions of $e\D(G)$ and $e\D_H(G)$} \label{altdes}
Let us now describe $e\D(G)$ and $e\D_H(G)$ as categories of certain quasi-equivariant complexes. The connected unipotent group $N$ acts by left translations on $G$. Let $\mu_N:N\times G\to G$ denote this action. This is the restriction of the multiplication map. Also, since $N$ is a normal subgroup of $H$, we have an action of $H$ on $N$ by conjugation, and we can form their semidirect product $H\ltimes N$. Then the action of $N$ on $G$ by left translations and the action of $H$ on $G$ by conjugation give us an action $\alpha$ of $H\ltimes N$ on $G$. Let $C:H\times G\to G$ denote the conjugation action. Since $\L$ is a $G$-equivariant local system on $N$, we see that $pr_2^*\L$ is in fact a multiplicative local system on $H\ltimes N$, where $pr_2:H\ltimes N\to N$ is the second projection. From now on, let us denote by $U$ the group $H\ltimes N$, and by $\L'$, the multiplicative local system $pr_2^*\L$ on $U$.
\bprop\label{quasieq}\mbox{}
\begin{itemize}
\item[(a)] $e\D(G)=\D_{N,\L}(G)$ as full subcategories of $\D(G)$.
\item[(b)] $e\D_H(G)=\D_{U, \L'}(G)$ as full subcategories of $\D(G)$.
\end{itemize}
\eprop
\bpf
(a) Let $M\in \D(G)$. Let us compute $\mu_N^*(e\ast M)$. Since $e$ is supported on $N$, we have that $$\mu_N^*(e\ast M)\cong \mu_N^*{\mu_N}_!(e\boxtimes M).$$ By proper base change, we have that $$\mu_N^*{\mu_N}_!(e\boxtimes M) \cong (id_N\times \mu_N)_!(\mu|_{N\times N}\times id_G)^*(\L[2\dim N](\dim N)\boxtimes M)$$ $$\cong (id_N\times \mu_N)_!(\L\boxtimes e\boxtimes M)$$ $$\cong \L\boxtimes (e\ast M).$$ 
Hence we have that $\mu_N^*(e\ast M)\cong \L\boxtimes (e\ast M)$, i.e. $e\ast M\in \D_{N,\L}(G)$.

On the other hand, if we have $\mu_N^*M\cong \L\boxtimes M$, then we have that $$\L\ast M\cong {\mu_N}_!\mu_N^*M\cong M\otimes{\mu_N}_!\Qlcl\cong M[-2\dim N](-\dim N).$$ The middle isomorphism is given by the projection formula. Hence we see that $e\ast M\cong M$. Hence we have that $e\D(G)=\D_{N,\L}(G)$.
\medbreak
\noin (b) We have that $e\D_H(G)= e\D(G)\cap\D_H(G)$. Hence, we see that $e\D_H(G)=\D_{N,\L}(G)\cap \D_{H,\Qlcl}(G)$. Hence it is clear that $\D_{U, \L'}(G)\subset e\D_H(G)$.

On the other hand, suppose $M\in \D_{N,\L}(G)\cap \D_{H,\Qlcl}(G)=e\D_H(G)$. The map $\alpha:(H\ltimes N)\times G\to G$ factors as $(H\ltimes N)\times G\stackrel{(pr_2,c)}{\longrightarrow} N\times G\stackrel{\mu_N}{\longrightarrow} G$, where the first map is given by $(pr_2,c):((h,n),g)\mapsto (n,hgh^{-1})$. Hence we see that $$\alpha^*M\cong (pr_2,c)^*\mu_N^*M\cong(pr_2,c)^*(\L\boxtimes M)$$ $$\cong (\Qlcl\boxtimes \L)\boxtimes M\cong \L'\boxtimes M.$$ Hence $M\in \D_{U,\L'}(G)$.
\epf

\section{The category $e\D_H(G)$}
In this section, we will study the $\Qlcl$-linear category $e\D_H(G)$ and give the proof of Theorem \ref{Main1}. We have seen that $e\D_H(G)=\D_{U,\L'}(G)$ as full subcategories of $\D(G)$. 

\subsection{Support of objects of $e\D_H(G)$} \label{supportedhg}
Let us show that there exist finitely many $U$-orbits in $G$, such that the support of every
object of $e\D_H(G)$ is contained in their union. First, for a $g\in G$, let us describe the stabilizer $U_g$. Note that we have $(h,n)\cdot g=nhgh^{-1}$. Hence $(h,n)\in U_g$ if and only if $hgh^{-1}g^{-1}=n^{-1}$. Let $c_g:H\to H$ be the commutator map defined by $c_g(h)=hgh^{-1}g^{-1}$. Note that we have the following identity
\beq\label{id}
c_{g_1g_2}(h)=c_{g_1}(h)\cdot {^{g_1}c_{g_2}}(h).
\eeq
Since $H/N$ is commutative we have that in fact $c_{h'}:H\to N$ for $h'\in H$. From these observations, we obtain
\bprop\label{Ug}
$U_g=\{\left(h,c_g(h)^{-1}\right)|c_g(h)\in N\}$. Let $H_g=pr_1(U_g)=\{h\in H|c_g(h)\in N\}$. Hence, we have the map $c_g:H_g\to N$. Then $c_g^*\L$ is a multiplicative local system on $H_g$. Moreover, the subgroup $H_g$ depends only on the coset $Hg$. 
\eprop

\bprop\label{support1}
Let $M\in \D_{U,\L'}(X)$. Suppose $g\in G$ is such that $M_g\neq 0$. Then $c_g^*\L\cong \Qlcl$ as local systems, or equivalently, the multiplicative local system $c_g^*\L|_{H^0_g}$ on $H^0_g$ is trivial.
\eprop
\bpf
By Proposition \ref{support}, we must have $\L'|_{U_g}\cong \Qlcl$ as local systems. But we have an isomorphism $pr_1:U_g\to H_g$, and under this isomorphism $\L'|_{U_g}$ corresponds to $c_g^*(\L^{-1})$.
\epf

From (\ref{id}), we see that we have $c_{h'g}^*\L\cong (c_{h'}^*\L)|_{H_g}\otimes c_{g}^*\L$ for $h'\in H$, where we consider $c_{h'g}, c_g$ as maps from $H_g\to N$. Let us now fix a $g\in G$, and find all $h'g\in Hg$ such that $c_{h'g}^*\L$ is trivial, or equivalently $(c_{h'}^*\L)|_{H_g}\cong c_g^*(\L^{-1})$.

We will now need the construction described in Appendix A.13 of \cite{B}. Note that we have a connected unipotent group $H$, with a connected normal subgroup $N$ such that $[H,H]\subset N$. We also have an $H$-equivariant multiplicative local system $\L$ on $N$. Then this construction gives us a map $\phi_{\L}:(H/N)_{perf}\to (H/N)_{perf}^*$, where $(H/N)_{perf}^*$ is the Serre dual of $(H/N)_{perf}$. Note that we have the map $c_{h'}:H\to N$ for $h'\in H$. The map $\phi_\L$ is induced by the map $h'\mapsto c_{h'}^*\L$. By our hypothesis, the map $\phi_\L$ is an isogeny. Let $H_g^0$ denote the identity component of $H_g$. Note that the inclusion $i:H_g^0/N\hookrightarrow H/N$ gives us the surjective map $i^*:(H/N)_{perf}^*\rightarrow (H_g^0/N)_{perf}^*$. Note that the isomorphism $(c_{h'}^*\L)|_{H_g}\cong c_g^*(\L^{-1})$ of local systems exists if and only if $(c_{h'}^*\L)|_{H^0_g}\cong c_g|_{H^0_g}^*(\L^{-1})$. Now $c_g|_{H^0_g}^*(\L^{-1})$ gives us an element of $s\in (H^0_g/N)_{perf}^*$. From this, we see that we have the following:
\bprop\label{finorbits}
The objects of $\D_{U,\L'}(Hg)$ can only be supported on those $h'g\in Hg$ such that $$i^*(\phi_\L(Nh'))=s.$$ This defines a closed subvariety of $Hg$ made up of finitely many $U$-orbits in $Hg$.
\eprop
\bpf
By what we have said above, it follows that an isomorphism $c_{h'g}^*\L\cong \Qlcl$ exists if and only if $(c_{h'}^*\L)|_{H^0_g}\cong c_g|_{H^0_g}^*(\L^{-1})$, i.e if and only if $i^*(\phi_\L(Nh'))=s$. The set of all such $h'g$ defines a closed subvariety of $Hg$ which is stable under the action of $U$, and has dimension equal to $\dim N+\dim(H/N)-\dim(H^0_g/N)$. Moreover, all $U$-orbits in $Hg$ are closed and have dimension equal to the number above. Hence we see that this subvariety must consist of finitely many $U$-orbits in $Hg$.
\epf

\subsection{Proof of Theorem \ref{Main1}} \label{secmain1}
We can now prove Theorem \ref{Main1}. Indeed, from Proposition \ref{finorbits} above, we see that objects of $e\D_H(G)=\D_{U,\L'}(G)$ are supported on finitely many $U$-orbits in $G$. Proposition \ref{transact} describes the categories of quasi-equivariant complexes and perverse sheaves supported on a single orbit. In particular, we see that $\M^{perv}_e$ must be a semisimple abelian category and that $e\D_H(G)$ must be its bounded derived category. From Proposition \ref{transact} it also clear that all the simple objects in $\M^{perv}_e$ must be suitably shifted local systems supported on a closed subset, namely a $U$-orbit in $G$.

\subsection{The category $e\D_H(H)$} \label{tcedhh}
Let us now study the braided monoidal category $e\D_H(H)$ with unit object $e$.
\blem
The $U$-orbits in $H$ are the cosets $Nh'$ of $N$ in $H$. Hence the objects of $e\D_H(H)$ are supported on the cosets $Nk$ of $N$ such that $Nk\in K_\L = ker(\phi_L:(\V)_{perf}\to (\V)_{perf}^*)$. 
\elem
\bpf
Note that for $h'\in H$, $H_{h'}=H$, since $c_{h'}(H)\subset N$. Then the first statement is clear, since $(h,n)\cdot h'=nhh'h^{-1}=n[h,h']h'\in Nh'$. The second statement follows from Proposition \ref{support1}.
\epf

\bprop\label{ek}
Let $k\in K:=\{k\in H|c_k^*\L\cong \Qlcl\}$. (We have $K/N=K_\L$.)  Let $e^k$ denote the right translate of $e$ by $k\in H$. Then $e^k\in e\D_H(H)$.
\eprop
\bpf
First, let us check that $e^k\in e\D(H)$. Indeed, by Proposition \ref{properties}, we have $e\ast e^k\cong (e\ast e)^k\cong e^k$. Note that for this, we do not require $k$ to lie in $K$. Let us now show that $e^k\in \D_H(H)$. Let $C:H\times H\to H$ denote the conjugation action and let $P_i:H\times H\to H$ denote the respective projections. We will construct an isomorphism $C^*e^k\cong P_2^*e^k$. From the $G$-equivariant structure on $\L$, we get an isomorphism $C^*e\cong P_2^*e$. Note that we have a commutative diagram 
$$\xymatrix{H\times H\ar[d]_{(id,r_{k^{-1}})} \ar[r]^{C} &H\\
H\times H\ar[d]_{(C,c_k\circ P_1)} &H\ar[u]_{r_k}\\
H\times N.\ar[ur]_{\mu}}$$
Hence we get a sequence of canonical isomorphisms 
$$C^*e^k{\cong}(id,r_{k^{-1}})^*(C,c_k\circ P_1)^*\mu^*{r_k}^*e^k$$ 
$$\cong (id,r_{k^{-1}})^*(C,c_k\circ P_1)^*\mu^*e$$
$$\cong (id,r_{k^{-1}})^*(C,c_k\circ P_1)^*(e\boxtimes \L)$$
$$\cong (id,r_{k^{-1}})^*(C^*e\otimes P_1^*{c_k}^*\L)$$
$$\cong (id,r_{k^{-1}})^*(P_2^*e\otimes P_1^*{c_k}^*\L)$$
$$\cong (id,r_{k^{-1}})^*({c_k}^*\L\boxtimes e)$$
$$\cong {c_k}^*\L\boxtimes e^k.$$
Note that by assumption, we have a trivialization $c_k^*\L\cong \Qlcl$. Hence we conclude that we have $C^*e^k\cong P_2^*e^k$. Hence indeed we have $e^k\in eD_H(H)$.
\epf

\bprop\label{edhh}
Let $k\in K.$ Then $e^k[-\dim N]$ is the unique irreducible perverse sheaf (up to isomorphism) in $e\D_H(H)$ supported on $Nk$. In particular, the isomorphism class of $e^k$ only depends on the coset $Nk$. Hence if we choose a set of coset representatives $k_i\in K$ of $K/N$, $e^{k_i}[-\dim N]$ are all the irreducible perverse sheaves in $e\D_H(H)$ (up to isomorphism).
\eprop
\bpf
Since $e$ is supported on $N$, $e^k$ is supported on $Nk$. We have seen above that $e^k\in e\D_H(H)$, hence $e^k\in \D_{U,\L'}(Nk)$. Note that the stabilizer of $k\in Nk$ in $U$ is $U_k\cong H_k= H$. (See Proposition \ref{Ug}.) Since $H$ is connected, by Proposition \ref{transact}, we see that $\M^{perv}_{U,\L'}(Nk)\cong  {Vec}$. In particular, $\D_{U,\L'}(Nk)$ has only one irreducible perverse sheaf up to isomorphism. The proposition now follows, since $e^k[-\dim N]$ is perverse.
\epf

\bprop\label{edhhmod}
Let $M\in e\D_H(G)$ and $k\in K$. Then $e^{k}\ast M\cong M\ast e^{k}\cong M^k$. In particular, if $M\in \M_e$ and $k\in K$ (and hence $e^k\in \M_e$), then $M\ast e^k\in \M_e$.
\eprop
\bpf
Since $e^k\in \D_H(H)$, we have that $e^{k}\ast M\cong M\ast e^{k}$. By Proposition \ref{properties}, we have $M\ast e^k\cong (M\ast e)^k\cong M^k$.
\epf

\brk
Let $\M_{e,1}$ denote the full subcategory of $\M_e$ consisting of complexes supported on $H$. Proposition \ref{edhhmod} above shows that $\M_{e,1}$ is closed under convolution and also gives us the `multiplication table' for $\M_{e,1}$. Indeed, by Proposition \ref{edhh}, the simple objects of $\M_{e,1}$ are given by the $e^k$ for $k\in K$. The proposition tells us that $e^{k_1}\ast e^{k_2}\cong e^{k_1k_2}$ for $k_1,k_2\in K$.
\erk

\bdefn
Let $\C$ be a semisimple tensor category over an algebraically closed field of characteristic zero, with simple unit object. We say that $\C$ is a pointed category if all simple objects in $\C$ are invertible.
\edefn

\bcor
The category $\M_{e,1}$ is a monoidal category with unit object $e$. It is a pointed category.
\ecor

\subsection{Twists in the category $e\D_H(H)$} \label{titcedhh}

Note that for $k\in K$, we have constructed in Proposition \ref{ek} an isomorphism $C^*e^k\stackrel{\cong}{\longrightarrow}P_2^*e^k$. Pulling back this isomorphism via the diagonal $\Delta:H\to H\times H$, we get an automorphism $\theta_k=\theta_{Nk}$ of $e^k$. We will call this the twist of $e^k$. Since $e^k$ is a simple object, this is a number in $\Qlcl^*$. Let us compute the twists.  Namely, we show that  these twists give a polarization of a certain non-degenerate symmetric bimultiplicative form $K_\L\times K_\L\to \Qlcl^*$. Note that the isogeny $\phi_\L:(H/N)_{perf} \to (H/N)_{perf}^*$ gives us a skew-symmetric bimultiplicative local system on $H/N\times H/N$. Hence as described in \cite[\S A.10]{B}, we get a non-degenerate symmetric pairing $B:K_\L\times K_\L\to \Qlcl^*$.

\bprop
The twists $\theta_{Nk}$ give us a quadratic form $\theta:K_\L\to \Qlcl^*$. This quadratic form is a polarization of the pairing $B$ above, namely we have $\theta_{Nk_1k_2}\theta_{Nk_1}^{-1}\theta_{Nk_2}^{-1}=B(Nk_1,Nk_2)$.
\eprop
\bpf
Note that by the proof of \ref{ek}, we see that for every $k\in H$ we have a canonical isomorphism $C^*e^k\cong P_1^*c_k^*\L\otimes P_2^*e^k$. Let $c:H\times H\to N$ denote the commutator map. Then as described in \cite{B} A.13, the map $\phi_\L$ is induced by the bimultiplicative local system $c^*\L$ on $H\times H$. We have a unique trivialization $\rho:(c^*\L)|_{H\times K}\stackrel{\cong}{\longrightarrow}\Qlcl$ of bimultiplicative torsors on $H\times K$. We get the trivialization $\rho_k:c_k^*\L\stackrel{\cong}{\longrightarrow}\Qlcl$ by pulling back $\rho$ by the map $H\to H\times K$ given by $h\mapsto (h,k)$. Hence we get the $H$-equivariant structure isomorphism $C^*e^k\cong P_1^*c_k^*\L\otimes P_2^*e^k \stackrel{\cong}{\longrightarrow}\Qlcl\otimes P_2^*e^k$. To compute the twist $\theta_k$, it is sufficient to compute its restriction to the stalk of $e^k$ at the point $k$. So let us consider the composition $\Delta_k:\{k\}\hookrightarrow H\stackrel{\Delta}{\rightarrow}H\times H$. Pulling back the isomorphisms above and in $\ref{ek}$ by $\Delta_k$, we get the following automorphism $$e(1)=e(1)\otimes\Qlcl\stackrel{\cong}{\longrightarrow}e(1)\otimes\L(1)\stackrel{id\otimes \rho(k,k)}{\longrightarrow}e(1)\otimes\Qlcl=e(1).$$
Note that the pullback of the bimultiplicative local system $c^*\L$ on $H\times H$ by the diagonal is the trivial multiplicative local system on $H$. This trivialization comes from the isomorphism $\L(1)\stackrel{\cong}{\longrightarrow}\Qlcl$. Hence using Lemma A.26 from \cite{B}, we conclude that $\theta$ is indeed a quadratic form that gives a polarization of $B$.
\epf

\section{Convolution of perverse sheaves}
In this section, we will prove Theorem \ref{Main2}. Let $M, L$ be irreducible perverse sheaves in $\M^{perv}_e$. We want to show that $M\ast L[\dim N]$ is also perverse. Let us first show, as a consequence of Artin's theorem, that $M\ast L\in {}^p\D^{\geq \dim N}(G)$.

\subsection{A consequence of Artin's theorem} \label{acoat}
The following result is essentially due to M. Artin:
\bthm\label{artin}
If $f:X\to Y$ is an affine morphism of separated schemes of finite type over $\F$, the functor $f_*:\D(X)\to \D(Y)$ takes ${}^{p}\D^{\leq 0}(X)$ into ${}^{p}\D^{\leq 0}(Y)$. Hence by Verdier duality, the functor $f_!:\D(X)\to \D(Y)$ takes ${}^{p}\D^{\geq 0}(X)$ into ${}^{p}\D^{\geq 0}(Y)$.
\ethm

\medbreak
\noin Let us use this theorem to prove the following:

\bprop\label{artin1}
Let $M, L$ be irreducible perverse sheaves in $\M^{perv}_e$. Then $M\ast L\in { }^p\D^{\geq \dim N}(G)$.
\eprop
\bpf
We have seen that $M, L$ are suitably shifted local systems supported on $U$-orbits in $G$. It follows that $M\boxtimes L$ is a perverse sheaf on $G\times G$. Consider the free action of $N$ on $G\times G$ given by $n\cdot(g_1,g_2)=(g_1n^{-1},ng_2)$. Then the multiplication map $\mu:G\times G\to G$, factors as $G\times G\stackrel{\pi}{\rightarrow} N\backslash (G\times G)\stackrel{\tilde{\mu}}{\rightarrow} G$. Now we know that $\pi^{*}$ induces an equivalence of categories $\D(N\backslash (G\times G))\cong \D_N(G\times G)$. Also for any $M'\in \D(N\backslash (G\times G))$, we have that $\pi_!(\pi^*M')\cong M'\otimes \pi_!\Qlcl=M'[-2\dim N](-\dim N)$ by the Projection formula. Hence we see that $\pi_!$ takes $N$-equivariant perverse sheaves on $G\times G$ to $^p\D^{\dim N}(N\backslash(G\times G))$. By Artin's Theorem, we see that $\tilde{\mu}$ takes $^p\D^{\dim N}(N\backslash(G\times G))$ into $^p\D^{\geq \dim N}(G))$. It follows that $M\ast L=\mu_!(M\boxtimes L)\in {}^p\D^{\geq \dim N}(G)$.
\epf

\subsection{Convolving with the dual} \label{cwtd}

For a simple object $M\in \M_e$, let us compute $M^{\vee}\ast M$ and show that it lies in $\M_e$. We will use this preliminary computation to show that $\M_e$ is closed under convolution.

\blem
For any $M\in e\D_H(G)$, let $K_M=\{k\in K|M^k\cong M\}$. Then  $K_M$ is a closed subgroup of $K$ which contains $N$.
\elem
\bpf
It is clear that $K_M$ is a subgroup. By Proposition \ref{edhhmod}, we have that $M^k\cong M\ast e^k$. Since the isomorphism class of $e^k$ only depends on the coset $Nk$, we see that the isomorphism class of $M^k$ depends only on the coset $Nk$. Hence $N\subset K_M$. Since $K/N$ is finite, it follows that $K_M$ is a closed subgroup of $K$ containing $N$. 
\epf

\bprop\label{dualastm}
Let $M\in \M_e$ be a simple object supported on the $U$-orbit of $g\in G$. Then $M^\vee\ast M\in \M_{e,1}$. In fact $M^{\vee}\ast M\cong \bigoplus\limits_{Nk\in K_M/N}e^k$.
\eprop
\bpf
Since $M\in \D_{U,\L'}(Hg)$ is simple, we see that $M^{\vee}\in \D_{U,\L'}(Hg^{-1})$ is also simple. Hence it follows that $M^{\vee}\ast M \in \D_{U,\L'}(H)=e\D_H(H)$.
By Proposition \ref{dualityinhecke}, we have that $\Hom(M^{\vee}[m],M^{\vee})=\Hom(M^{\vee}\ast M[m],e)$ for any $m\in \f{Z}$. In other words, since $M^{\vee}$ is a simple object, we have that $\Hom(M^{\vee}\ast M,e)=\Qlcl$ and $\Hom(M^{\vee}\ast M[m],e)=0$ for $m\neq 0$. Also for any $k\in K$ and $m\in \f{Z}$, we have that $$\Hom(M^{\vee}\ast M[m],e^{k^{-1}})=\Hom(M^{\vee}\ast M\ast e^k[m], e)=\Hom(M^{\vee}\ast M^k[m],e)$$ $$=\Hom(M^{\vee}[m],(M^k)^{\vee}).$$ We see that
\[ \Hom(M^{\vee}[m],(M^k)^{\vee}) = \left\{ \begin{array}{ll}
         \Qlcl & \mbox{if } m=0 \mbox{ and } M^k\cong M\\
        0 & \mbox{otherwise}.\end{array} \right. \] 
Then using the fact that $\M_e$ is semisimple and that $e\D_H(G)$ is its bounded derived category, we conclude that we must have $$M^{\vee}\ast M\cong \bigoplus\limits_{Nk\in K_M/N}e^k.$$
In particular $M^{\vee}\ast M\in \M_{e,1}$.
\epf

\blem
Let $M,L\in \M_e$ be nonzero. Then $M\ast L$ is also nonzero. 
\elem
\bpf
We may assume that $M,L$ are simple objects. Suppose  $M\ast L= 0$. Then $(M^{\vee}\ast M)\ast L= 0$, i.e. $\left(\bigoplus\limits_{Nk\in K_M/N}e^k\right)\ast L= 0.$ By Proposition \ref{edhhmod}, we see that this is absurd.
\epf

\subsection{Proof of Theorem \ref{Main2}} \label{potm2}
Let us now complete the proof of Theorem \ref{Main2}. Let $M,L\in \M_e$ be simple. Then by Proposition \ref{artin1}, we know that $M\ast L\in {}^p\D^{\geq -\dim N}(G)$. Also, $M\ast L\in e\D_H(G)$ which is the bounded derived category of the semisimple abelian category $\M_e$. Hence we see that we must have $M\ast L\cong P^0\oplus P^1[-1]\cdots \oplus P^m[-m]$ for some non-negative integer $m$ and $P^i\in \M_e$. Now we have that $$M^{\vee}\ast (M\ast L)\cong M^{\vee}\ast P^0\oplus M^{\vee}\ast P^1[-1]\cdots \oplus M^{\vee}\ast P^m[-m].$$
On the other hand $(M^{\vee}\ast M)\ast L\in \M_e$ by Propositions \ref{dualastm} and \ref{edhhmod}. By Proposition \ref{artin1}, we have that $M^{\vee}\ast P^i[-i]\in {}^p\D^{\geq -\dim N+i}(G)$. Hence for $i>0$, we must have that $M^{\vee}\ast P^i\cong 0$. From the lemma, we conclude that we must have $P^i\cong 0$ for all $i>0$, i.e. $M\ast L\cong P^0$. In other words $M\ast L\in \M_e$. Also we have that $e\in \M_e$. Hence $\M_e$ is a full subcategory of the monoidal category $e\D_H(G)$ that is closed under convolution, and contains the unit $e$. Hence $\M_e$ is indeed a monoidal category with unit object $e$.

\section{Rigidity of $\M_e$} \label{rigidity}
In this section, we will prove that the category $\M_e$ is rigid. $\M_e$ is graded by the finite group $\Gamma$. This grading is faithful (Defn. \ref{grading}). Moreover, we have seen at the end of \S\ref{tcedhh} that the identity component $\M_{e,1}$ is pointed. Also, we have described a weak notion of duality in the category $\M_e$. Hence the rigidity of $\M_e$ follows from Theorem \ref{rigid} below.

\subsection{Rigidity in certain graded tensor categories} \label{ricgtc}
Let $\C$ be a tensor category over a field $k$ of characteristic zero such that:
\bit
	\item[(i)] As a $k$-linear category, $\C\cong Vec\oplus\cdots \oplus Vec.$
	\item[(ii)] $\mbox{End} \mathbf{1}=k$.
	\item[(iii)] For every simple object $M\in \C$, there exists a simple object $M^\vee \in \C$ such that\footnote{Under the first two conditions, this property is equivalent to the weak notion of duality described in \S\ref{catdg}.} $$\dim\Hom(M\otimes M^\vee, \mathbf{1})=\dim\Hom(M^\vee\otimes M, \mathbf{1})=1$$ and $$\Hom(M\otimes Y, \mathbf{1})=\Hom(Y\otimes M, \mathbf{1})=0$$ for all simple objects $Y\in \C$ not isomorphic to $M^\vee$.
	\item[(iv)] $\C$ has a grading $\C=\bigoplus\limits_{\gamma\in \Gamma}{\C_\gamma}$ by a finite group $\Gamma$ so that $\C_1$ is pointed. Let $\mathcal{G}$ denote the group of isomorphism classes of simple objects of $\C_1$.
\eit
\bthm\label{rigid}
Let $\C$ be a tensor category as above. Then $\C$ is rigid and hence a fusion category.
\ethm

\noin We will prove this theorem in \S\ref{potr}. Let us first recall some facts about rigidity and duals.

\subsection{Rigidity}
\bdefn
Let $\C$ be a monoidal category. Let $M$ be an object of $\C$. A left dual of $M$ is a triple $(M^*,ev_M,coev_M)$, where $M^*$ is an object of $\C$, $ev_M : M^*\otimes M\to \mathbf{1}$ and $coev_M:\mathbf{1}\to M\otimes M^*$ such that the compositions 
\beq \label{comp1}
M\cong \mathbf{1}\otimes M\stackrel{coev_M\otimes id_M}{\longrightarrow}(M\otimes M^*)\otimes M\cong M\otimes(M^*\otimes M)\stackrel{id_M\otimes ev_M}{\longrightarrow}M\otimes \mathbf{1}\cong M
\eeq
 and 
\beq\label{comp2}
M^*\cong M^*\otimes\mathbf{1}\stackrel{id_{M^*}\otimes coev_M}{\longrightarrow}M^*\otimes(M\otimes M^*)\cong (M^*\otimes M)\otimes M^*\stackrel{ev_M\otimes id_{M^*}}{\longrightarrow} \mathbf{1}\otimes M^* \cong M^*
\eeq
 are equal to the identity morphisms.
\edefn

\brk\label{leftadj}
\mbox{ }Given a triple $(M^*,ev_M,coev_M)$, we can define \mbox{ }maps in either direction \mbox{ }between $\Hom(A, M \otimes B)$ and $\Hom(M^*\otimes A, B)$ for all $A,B\in \C$. The conditions (\ref{comp1}) and (\ref{comp2}) imply that the compositions of these in either direction equal the identity maps. Hence in this case, the functor $M^*\otimes(\cdot)$ is left adjoint to $M\otimes (\cdot)$.
\erk

\blem\label{identity}
Let $\C$, $M$ be as above. Let $(M^*,ev_M,coev_M)$ be a triple such that one of the compositions (\ref{comp1}) and (\ref{comp2}) is the identity, while the other one is an isomorphism. Then the other composition must also be the identity.
\elem
\bpf
Let $A$, $B\in \C$. Then, as above, we get maps in either direction between  $\Hom(A, M \otimes B)$ and $\Hom(M^*\otimes A, B)$.  Since one of (\ref{comp1}) and (\ref{comp2}) is the identity, while the other one is an isomorphism, we can deduce that composition of these maps in one direction is the identity, while it is an isomorphism in the other direction. Hence it follows that the composition in the reverse order must also be identity. From this we can deduce that both (\ref{comp1}) and (\ref{comp2}) must be identity morphisms.
\epf

\subsection{Proof of Theorem \ref{rigid}} \label{potr}
In this section, we will prove Theorem \ref{rigid}. I would like to thank M. Boyarchenko for simplifying the proof. Let $\mathcal{S}$ be a representative system of isomorphism classes of simple objects of $\C$ such that $\mathbf{1}\in \mathcal{S}$. We will consider $\mathcal{G}$ as a subset of $\mathcal{S}$.  We have a grading $\C=\bigoplus\limits_{\gamma\in \Gamma}\C_\gamma$. Let us fix some $M\in \cS$ lying in $\C_\gamma$. Then we must have $M^\vee\in \C_{\gamma^{-1}}$. To prove rigidity of $\C$, it is enough to prove that $M^\vee$ is in fact the (left) dual of $M$. We may assume that $M^\vee\in \cS$. Let us now compute the objects $M\otimes M^\vee, M^\vee\otimes M\in \C_1$.

\blem
\bit
\item[(i)] $X\in \cG$ occurs in $M\otimes M^\vee$ iff $X\otimes M\cong M$ iff $M^\vee\otimes X\cong M^\vee$. In this case, it occurs with multiplicity one.
\item[(ii)] $Y\in \cG$ occurs in $M^\vee\otimes M$ iff $M\otimes Y\cong M$ iff $Y\otimes M^\vee\cong M^\vee$. In this case, it occurs with multiplicity one.
\eit
\elem
\bpf
For $X\in \cG$, we have $\Hom(M\otimes M^\vee, X)\cong \Hom(X^{-1}\otimes M \otimes M^\vee, \mathbf{1})$ which is 0 unless $X^{-1}\otimes M\cong M\Longleftrightarrow X\otimes M\cong M$ in which case it is a $1$-dimensional vector space, i.e. in this case $X$ occurs with multiplicity one in $M\otimes M^\vee$. The other assertions are similar.
\epf

\noin Let $\A=\{X\in \cG|X\otimes M\cong M\}$ and $\B=\{Y\in \cG|M\otimes Y\cong M\}$. For each $X\in \A$, we fix a nonzero map (which is unique up to scaling) $c_X:X\to M\otimes M^\vee$.  Similarly, for each $Y\in \B$ fix a nonzero map $e_Y : M^\vee\otimes M{\longrightarrow} Y$. 

\blem\label{easier}
There exists $X\in \A$ such that the composition \beq\label{comp3} X\otimes M\stackrel{c_X\otimes id_M}{\longrightarrow} M\otimes M^\vee \otimes M\stackrel{id_M\otimes e_{\mathbf{1}}}{\longrightarrow} M\otimes {\mathbf{1}}\eeq is nonzero.
\elem
\bpf
The map $M\otimes M^\vee \otimes M\stackrel{id_M\otimes e_{\mathbf{1}}}{\longrightarrow} M\otimes {\mathbf{1}}$ is nonzero. Also, we have $M\otimes M^\vee \otimes M\cong \bigoplus\limits_{X\in \A}X\otimes M$ with the inclusions being given by $c_X\otimes id_M$. Hence we conclude that for some $X\in \A$ the composition (\ref{comp3}) must be nonzero.
\epf
\noin In fact, let us now prove that we can take $X=\mathbf{1}$.


\bprop
For the triple $(M^\vee,e_{\mathbf{1}},c_{\mathbf{1}})$, the compositions (\ref{comp1}) and (\ref{comp2}) are isomorphisms.
\eprop
\bpf
Let us prove that (\ref{comp1}) is an isomorphism. A similar argument can be used to prove that (\ref{comp2}) is an isomorphism.
Since $M$ is a simple object, it is enough to show that (\ref{comp1}) is nonzero, i.e. the composition $$\mathbf{1}\otimes M\stackrel{c_{\mathbf{1}}\otimes id_M}{\longrightarrow} M\otimes M^\vee \otimes M\stackrel{id_M\otimes e_{\mathbf{1}}}{\longrightarrow} M\otimes {\mathbf{1}}$$ is nonzero. Let $X\in \A$ be as in Lemma \ref{easier}. Let us fix an isomorphism $\phi:X^{-1}\otimes M \stackrel{\cong}{\longrightarrow} M$. We have the following commutative diagram:
\beq\label{diagram1}\xymatrix{X^{-1}\otimes M\otimes M^\vee\otimes M\ar[d]^{\cong}_{\phi\otimes id_{M^\vee}\otimes id_M} \ar[r]^{\mbox{  }\mbox{  }\mbox{  }\mbox{  }\mbox{  }id \otimes e_{\mathbf{1}}} &X^{-1}\otimes M\otimes \mathbf{1}\ar[d]_{\cong}^{\phi\otimes id_{\mathbf{1}}}\\
M\otimes M^\vee \otimes M \ar[r]_{\mbox{ }\mbox{ }\mbox{ }id_M\otimes e_{\mathbf{1}}} & M\otimes \mathbf{1} .}\eeq Also, since $\dim\Hom(\mathbf{1},M\otimes M^\vee)=1$, there exists an isomorphism $\psi:X^{-1}\otimes X \stackrel{\cong}{\longrightarrow} \mathbf{1}$ such that the following diagram commutes:

\beq\xymatrix{X^{-1}\otimes X\ar[d]^{\cong}_{\psi} \ar[r]^{id \otimes c_X\mbox{ }\mbox{ }\mbox{ }\mbox{ }} &X^{-1}\otimes M\otimes M^\vee\ar[d]_{\cong}^{\phi\otimes id_{M^\vee}}\\
\mathbf{1} \ar[r]_{c_{\mathbf{1}}\mbox{ }\mbox{ }\mbox{ }\mbox{ }} & M\otimes M^\vee .}\eeq Tensoring this diagram by $M$ and putting it together with (\ref{diagram1}), we obtain:
\beq\xymatrix{X^{-1}\otimes X\otimes M\ar[d]^{\cong}_{\psi\otimes id_M} \ar[r]^{id \otimes c_X\otimes id\mbox{ }\mbox{ }\mbox{ }\mbox{ }} &X^{-1}\otimes M\otimes M^\vee\otimes M\ar[d]_{\cong}^{\phi\otimes id_{M^\vee\otimes M}}\ar[r]^{\mbox{  }\mbox{  }\mbox{  }\mbox{  }\mbox{  }id \otimes e_{\mathbf{1}}} &X^{-1}\otimes M\otimes \mathbf{1}\ar[d]_{\cong}^{\phi\otimes id_{\mathbf{1}}}\\
\mathbf{1}\otimes M \ar[r]_{c_{\mathbf{1}}\otimes id_M\mbox{ }\mbox{ }\mbox{ }\mbox{ }\mbox{ }} & M\otimes M^\vee\otimes M \ar[r]_{\mbox{ }\mbox{ }\mbox{ }id_M\otimes e_{\mathbf{1}}} & M\otimes \mathbf{1}.}\eeq
By Lemma \ref{easier}, the composition of the two maps on the top row of this diagram is an isomorphism. Hence the composition in the bottom row must be an isomorphism. Hence we see that (\ref{comp1}) is an isomorphism.

\epf

\bpf[Proof of Theorem \ref{rigid}.]
Since (\ref{comp1}) is nonzero, we can rescale the map $c_{\mathbf{1}}$ and ensure that (\ref{comp1}) is the identity. Then since (\ref{comp2}) is an isomorphism, by Lemma \ref{identity} it must also be the identity. Hence $(M^\vee,e_{\mathbf{1}},c_{\mathbf{1}})$ is dual to $M$. Hence all $M\in \cS$ have a dual. Hence $\C$ is a fusion category.
\epf

\noin Hence we conclude that $\M_e$ is a fusion category, completing the proof of Theorem \ref{Main3}.

\section{Modularity of $\M_e^{\Gamma}$} \label{modularity}
As we have noted earlier, we have a twist automorphism $\theta$ of the identity functor on $\D_G(G)$. $\D_G(G)$ is a braided monoidal category. Let $\beta$ denote the braiding. Let $(M,\phi)$ be an object of $\D_G(G)$. Then $\theta_{(M,\phi)}:=\Delta^*\phi:(M,\phi)\to (M,\phi)$, where $\Delta:G\to G\times G$ is the diagonal. (See \cite[\S 3.9]{B}.) These twists satisfy the following {\it balancing property} 
\beq
\theta_{M\ast L}=\beta_{L,M}\circ \beta_{M,L}\circ (\theta_M\ast \theta_L), \mbox{ for all } M,L\in \D_G(G).
\eeq
In this section we will show that $\Meg$ with the twist $\theta$ defined above is a modular category. First let us show that it is a ribbon category.

\subsection{Ribbon property} \label{ribbon}
Let us recall the definition of a ribbon structure.
\bdefn\label{ribbondefn}
Let $\C$ be a rigid braided monoidal category with braiding $\beta$. A ribbon structure on $\C$ is an automorphism $\theta$ of the identity functor on $\C$ satisfying the following two conditions:
\bit
\item[(i)] $\theta_{M\otimes L}=\beta_{L,M}\circ \beta_{M,L}\circ (\theta_M\ot \theta_L), \mbox{ for all } M,L\in \C$,
\item[(ii)] $\theta_{M^{\ast}}=\theta_M^{\ast}$ for all $M\in \C$.
\eit
\edefn

\noin Let us now show that the twist $\theta$ defined on $\Meg$ defines a ribbon structure on $\Meg$. As before, let $\iota:G\to G$ denote the inversion map. 
\bprop
The twists satisfy the following relations:
\bit
\item[(i)] $\iota^*\theta_{(M,\phi)}={(\theta_{(\iota^*M, (id_G\times \iota)^*\phi)})}^{-1}$ for all $(M,\phi)\in\D_G(G)$.
\item[(ii)] $\mathbb{D}\theta_{(M,\phi)}={(\theta_{(\mathbb{D}M,\mathbb{D}\phi^{-1}[2\dim G](\dim G))})}^{-1}$ for all $(M,\phi)\in\D_G(G)$.
\eit
Hence $\mathbb{D}^-\theta_M = \theta_{\mathbb{D}^-M}$ for all $M\in \D_G(G)$ and $\theta_M^\vee=\theta_{M^\vee}$ for all $M\in e\D_G(G)$. Since $\Meg\subset e\D_G(G)$ is a fusion category with duality defined by $(\cdot)^\vee$, we conclude that $\theta$ defines a ribbon structure on $\Meg$.
\eprop
\bpf
(i) follows from the equality $\iota^*\Delta^*\phi = \Delta^*(id_G\times \iota)^*\phi^{-1}$, which is a result of the compatibility relation satisfied by $\phi$. To prove (ii), first note that $\theta_{M\otimes L}=\theta_M\otimes \theta_L$ for $M,L\in \D_G(G)$. Also, we have functorial isomorphisms $\Hom(L,\mathbb{D}M)\cong \Hom(L\otimes M, \K)$ for $M,L\in \D_G(G)$, where $\K=\K_G$ is the dualizing sheaf on $G$. Under this correspondence, let $ev:\f{D}M\otimes M\to \K$ be the map corresponding to $id_{\f{D}M}$. Then the isomorphism $\Hom(\f{D}M,\f{D}M)\cong \Hom(\f{D}M\otimes M, \K)$ is given by $f \mapsto ev\circ (f\otimes id_M)$. Since $\theta_\K=id_\K$ and since $\theta$ is an automorphism of the identity functor, we have that $ev\circ(\theta_{\f{D}M\otimes M})=\theta_\K\circ ev$ and hence $ev\circ (id_{\f{D}M}\otimes \theta_M) = ev\circ (\theta_{\f{D}M}^{-1}\otimes id_M)$. Hence we conclude that indeed $\f{D}\theta_M=\theta_{\f{D}M}^{-1}$.
Since $\f{D}^-=\iota^*\circ\f{D}=\f{D}\circ\iota^*$, we conclude from (i) and (ii) that $\mathbb{D}^-\theta_M = \theta_{\mathbb{D}^-M}$ for all $M\in \D_G(G)$. And since $(\cdot)^\vee=\f{D}^{-}(\cdot)[2\dim N](\dim N)$ we also have $\theta_M^\vee=\theta_{M^\vee}$ for all $M\in e\D_G(G)$.
\epf

\subsection{Proof of Theorem \ref{Main4}} \label{nondeg}
Theorem \ref{Main4}(i) follows from Theorem \ref{Main1}. We have seen that the category $\M_e$ is a braided $\Gamma$-crossed category. Hence it follows that $\Meg$ has the structure of braided monoidal category. Rigidity of $\Meg$ follows from the rigidity of $\M_e$. Combining these observations with Theorem \ref{Main4}(i), we deduce Theorem \ref{Main4}(ii). Let us now prove statement (iii). It follows from \S\ref{ribbon} that the identity component $\M_{e,1}$ is a ribbon category. In \S\ref{titcedhh}, we have seen that the twist on $\M_{e,1}$ is given  by a quadratic form which gives a polarization of the non-degenerate symmetric pairing $B:K_\L\times K_\L\to \Qlcl^*$. Since $B$ is non-degenerate, it follows that $\M_{e,1}$ is a modular category. Then it follows from \cite[Prop. 4.56(ii)]{DGNO} that $\Meg$ must be a non-degenerate braided category. We have seen above that $\Meg$ is a pre-modular category. Hence we conclude that $\Meg$ is a modular category.

\bibliographystyle{ams-alpha}

\end{document}